\documentclass[11pt]{article}%
\usepackage{amsmath}
\usepackage{amsfonts}
\usepackage{amssymb}
\usepackage{graphicx}
\usepackage{a4wide}%
\providecommand{\U}[1]{\protect\rule{.1in}{.1in}}

\begin{document}

\title{On Berinde's paper\\ \textquotedblleft Comments on some fixed point
theorems in metric spaces"}
\author{C. Z\u{a}linescu\thanks{Octav Mayer Institute of Mathematics, Ia\c{s}i Branch
of Romanian Academy, Ia\c{s}i, Romania, email: \texttt{zalinesc@uaic.ro.}}}
\date{}
\maketitle

\begin{abstract}
In Theorem 1 of the paper [V. Pata, A fixed point theorem in metric
spaces, J. Fixed Point Theory Appl., 10 (2011), 299--305] it is
proved that Picard's iterates for a function converge to a fixed
point if a certain condition (C) is verified for all parameters in
the interval [0,1]. In the recent paper [V. Berinde, Comments on
some fixed point theorems in metric spaces, Creat. Math. Inform. 27
(2018), 15--20], the author claims that Pata's result does not hold
at least in the two extremal cases for the parameter involved in
(C). In this note we point out that Berinde's Theorem 1.1 has only a
visual similarity with Pata's Theorem 1, and the conclusion of
Theorem 1.1 is verified only by constant functions.
\end{abstract}

The paper mentioned in the title, that is \cite{Ber18}, begins with the
following text:

\smallskip \textquotedblleft Let $(X,d)$ be a metric space. By selecting an
arbitrary point $x_{0}\in X$, which we call the zero of the metric
space $X$, we denote, according to the terminology and notations in
[19],

$\left\Vert x\right\Vert :=d(x,x_{0})$ $\forall x\in X$. \noindent
We also consider an increasing function $\Psi:[0,1]\rightarrow
\lbrack0,\infty)$ which is vanishing with continuity at zero and the
(vanishing) sequence

$\omega_{n}(\alpha)=(\frac{\alpha}{n})\sum_{k=1}^{n}\Psi(\frac{\alpha}{k})$,
\ \ $(1.1)$

\noindent where $\alpha\geq1$.

The following theorem is the main result in [19].

\smallskip {\bf Theorem 1.1.} {\it Let $(X,d)$ be a complete metric space and
$f:X\rightarrow X$ a self mapping of $X$. Let $\Lambda\geq0$,
$\alpha\geq1$ and $\beta\in\lbrack0,\alpha]$ be fixed constants. If
the inequality

$d(f(x),f(y))\leq(1-\varepsilon)d(x,y)+\Lambda\varepsilon^{\alpha}%
\Psi(\varepsilon)[1+\left\Vert x\right\Vert +\left\Vert y\right\Vert
]\ \ \ (1.2)$

\noindent is satisfied for every $\varepsilon\in\lbrack0,1]$ and
every $x,y\in X$, then $f$ possesses a unique fixed point
$x^{\ast}=f(x^{\ast})$. Furthermore, by denoting the $n$th iterate
of $f$ by $f^{n}$, we have the estimate

$d(x^{\ast},f^{n}(x_{0}))\leq C\omega_{n}(\alpha), \ \ (1.3)$

\noindent for some positive constant $C\leq\Lambda(1+4\left\Vert
x^{\ast }\right\Vert )\beta$.}

\smallskip Our aim in this note is to show that Theorem 1.1 does not hold at
least for two extremal cases of the parameter $\varepsilon$ involved in the
contraction condition (1.2). We also provide a correct version (but not fully
in the spirit) of Theorem 1.1 and discuss some other related results."

\medskip The reference \textquotedblleft\lbrack19]" mentioned in the quoted
text is Pata's paper \cite{Pat11}. The result envisaged by Berinde
in \cite{Ber18} is the following theorem quoted from \cite[page
299]{Pat11}:

\smallskip\textquotedblleft {\bf Theorem 1.} {\it Let $\Lambda\geq0$, $\alpha\geq1$ and
$\beta\in\lbrack0,\alpha]$ be fixed constants. If the inequality

$d(f(x),f(y))\leq(1-\varepsilon)d(x,y)+\Lambda\varepsilon^{\alpha}%
\psi(\varepsilon)[1+\left\Vert x\right\Vert +\left\Vert y\right\Vert
]^{\beta }\ \ \ (1.1)$

\noindent is satisfied for every $\varepsilon\in\lbrack0,1]$ and every $x,y\in
X$, then $f$ possesses a unique fixed point $x_{\ast}=f(x_{\ast})$.
Furthermore, calling $f^{n}=f\circ\cdots\circ f$ ($n$ times),

$d(x_{\ast},f^{n}(x_{0}))\leq C\omega_{n}(\alpha)\ \ \ (1.2)$

\noindent for some positive constant $C\leq\Lambda(1+4\left\Vert
x_{\ast }\right\Vert )^{\beta}$.}"

\smallskip In \cite[Theorem 1]{Pat11} $(X,d)$, $f$ and $\psi$ are as
in \cite[Theorem 1.1]{Ber18}, while
$\omega_{n}(\alpha)$ is $(\frac{\alpha}{n})^{\alpha}\sum_{\kappa=1}^{n}%
\psi(\frac{\alpha}{\kappa})$.

\medskip Comparing the statements of \cite[Theorem 1]{Pat11} and \cite[Theorem
1.1]{Ber18} we emphasize the following differences:

\smallskip(a)$~[1+\left\Vert x\right\Vert +\left\Vert y\right\Vert ]^{\beta}$
in \cite{Pat11} versus $[1+\left\Vert x\right\Vert +\left\Vert y\right\Vert ]$
in \cite{Ber18};

\smallskip(b)$~C\leq\Lambda(1+4\left\Vert x_{\ast}\right\Vert )^{\beta}$ in
\cite{Pat11} versus $C\leq\Lambda(1+4\left\Vert x^{\ast}\right\Vert )\beta$ in
\cite{Ber18};

\smallskip(c)$~\omega_{n}(\alpha)=(\frac{\alpha}{n})^{\alpha}\sum_{\kappa
=1}^{n}\psi(\frac{\alpha}{\kappa})$ in \cite{Pat11} versus $\omega_{n}%
(\alpha)=(\frac{\alpha}{n})\sum_{\kappa=1}^{n}\Psi(\frac{\alpha}{k})$ in
\cite{Ber18}.

\medskip Of course, one could (reasonably) assume that the differences (a),
(b) and (c) mentioned above are misprints.

\smallskip The fact that the missing $\beta$ in \cite[Eq.\ (1.2)]{Ber18} is
not a misprint is proved by the following text from \cite[page 17]{Ber18}:

\smallskip\textquotedblleft If in (1.2) we have $\varepsilon=1$ (or we let
$\varepsilon\rightarrow1$), then this condition becomes

$d(f(x),f(y))\leq L\cdot\lbrack1+d(x,x_{0})+d(y,x_{0})],\ \forall x,y\in
X,\ \ \ \ $(3.10)

\noindent for a fixed element $x_{0}\in X$ and a constant $L=\Lambda\psi(1)$".

\medskip

In the sequel we denote $\Psi$ from \cite{Ber18} with $\psi$. Having
in view the above differences, the natural question is if
\cite[Theorem 1]{Pat11} and \cite[Theorem 1.1]{Ber18} are
equivalent.

\smallskip A first remark is that in \cite[Theorem 1]{Pat11} and \cite[Theorem
1.1]{Ber18} $\omega_{n}(\alpha)$ does not make sense if $\alpha>1$
because $\psi(\frac{\alpha}{\kappa})$ is not well defined at least
for $\kappa=1$. The solution is simple: set $\psi(t):=\psi(1)$ for
$t>1$, or to take from the beginning
$\psi:\mathbb{R}_{+}:=[0,\infty\lbrack{}\rightarrow \mathbb{R}_{+}$
with the mentioned properties; this is practically done in
\cite[Corollary 2]{Pat11} and its proof by taking
$\psi(\varepsilon)=\varepsilon^{\gamma}$ (with $\gamma>0$ and,
without mentioning it explicitly, $\varepsilon\geq0$).

\smallskip We did not detect other drawbacks in the proof of \cite[Theorem
1]{Pat11}.

\medskip For further discussions, we rewrite the condition that $(X,f,\alpha
,\psi,\Lambda,\beta)$ (with $(X,d)$ metric space, $f:X\rightarrow
X$, $\alpha\geq1$, $\psi:\mathbb{R}_{+}\rightarrow \mathbb{R}_{+}$
increasing and continuous at $0$ with $\psi(0)=0$, $\Lambda\geq0$,
$\beta\in\lbrack0,\alpha]$) from
\cite[Theorem 1.1]{Ber18} must satisfy in the form%
\begin{equation}
\forall\varepsilon\in\lbrack0,1],\ \forall x,y\in X:d(f(x),f(y))\leq
(1-\varepsilon)d(x,y)+\Lambda\varepsilon^{\alpha}\psi(\varepsilon
)[1+\left\Vert x\right\Vert +\left\Vert y\right\Vert ].\label{r-b1.2}%
\end{equation}

So, $(X,f,\alpha,\psi,\Lambda,\beta)$ satisfies (\ref{r-b1.2}) if and only if
$(X,f,\alpha,\psi,\Lambda,\beta=0)$ satisfies (\ref{r-b1.2}). Having in view
this remark we get the following statement.

\medskip {\bf Fact B.} {\it Assume that \cite[Theorem 1.1]{Ber18} is true. If
$(X,f,\alpha,\psi,\Lambda,\beta)$ verifies the hypothesis of
\cite[Theorem 1.1]{Ber18}, then $f$ is constant.}

\medskip Proof. As observed above, $(X,f,\alpha,\psi,\Lambda,\beta=0)$
satisfies (\ref{r-b1.2}). By \cite[Theorem 1.1]{Ber18} there exists
a unique $\overline{x}\in X$ such that
$f(\overline{x})=\overline{x}$ and $d(\overline{x},f^{n}(x_{0}))\leq
C\omega_{n}(\alpha)$ for every $n\geq1$, where $0\leq
C\leq\Lambda(1+4\left\Vert x^{\ast}\right\Vert )\beta=0$. Therefore,
$f(x_{0})=\overline{x}$. Fix now $x_{0}^{\prime}\in X$ and set
$\left\Vert x\right\Vert ^{\prime }:=d(x,x_{0}^{\prime})$. Then for
$x,y\in X$ one has
\begin{align*}
1+\left\Vert x\right\Vert +\left\Vert y\right\Vert  &  =1+d(x,x_{0}%
)+d(y,x_{0})\leq1+2d(x_{0},x_0^{\prime})+d(x,x_0^{\prime})+d(y,x_0^{\prime})\\
&  \leq\left[  1+2d(x_{0},x_0^{\prime})\right]
\cdot\lbrack1+\left\Vert x\right\Vert ^{\prime}+\left\Vert
y\right\Vert ^{\prime}].
\end{align*}
It follows that \cite[Eq.\ (1.2)]{Ber18} holds for $x_{0}$,
$\Lambda$, and $\beta=0$ replaced by $x_{0}^{\prime}$,
$\Lambda^{\prime}:=\Lambda\cdot \lbrack1+2d(x_{0},x_0^{\prime})]$,
$\beta^{\prime}:=0$ (and the same $\alpha$ and $\psi$). Applying
again \cite[Theorem 1.1]{Ber18}, we get a (unique)
$\overline{x}^{\prime}\in X$ such that
$f(\overline{x}^{\prime})=\overline {x}^{\prime}=f(x_{0}^{\prime})$.
By the uniqueness of the fixed point $\overline{x}$ we obtain that
$f(x_{0}^{\prime})=\overline{x}$. Since $x_{0}^{\prime}$ is
arbitrary, it follows that $f$ is constant.

\medskip Note that any $\lambda$-contraction $f$ verifies condition
(\ref{r-b1.2}) for $\alpha:=1$, $\psi(\varepsilon):=\varepsilon$ for
$\varepsilon\geq0$, $\Lambda:=[4(1-\lambda)]^{-1}$ and $\beta:=0$;
indeed $1-\varepsilon+\Lambda
\varepsilon^{2}\geq1-\frac{1}{4\Lambda}=\lambda$ for $\varepsilon\in
\mathbb{R}$, and so
\[
d(f(x),f(y))   \leq\lambda d(x,y)\leq\left[  1-\varepsilon+\Lambda
\varepsilon^2\right]  d(x,y)
\leq(1-\varepsilon)d(x,y)+\Lambda\varepsilon^{\alpha}\psi(\varepsilon
)[1+\left\Vert x\right\Vert +\left\Vert y\right\Vert ]\] for all
$\varepsilon\in\lbrack0,1]$ and all $x,y\in X$. So, we conclude that
\cite[Theorem \ 1.1]{Ber18} is false.

\medskip Therefore, \cite[Theorem 1]{Pat11} and \cite[Theorem 1.1]{Ber18} are
{\bf not} equivalent.

\medskip Recall that Berinde's aim in his \textquotedblleft note is to show
that Theorem 1.1 does not hold at least for two extremal cases of the
parameter $\varepsilon$ involved in the contraction condition (1.2)," aim
which was achieved as mentioned in \cite[page 18]{Ber18}: \textquotedblleft By
summarizing the comments above, we conclude that Theorem 1.1 does not hold if
we have $\varepsilon=0$ or $\varepsilon=1$ in (1.2)".

\medskip Indeed, Berinde proved that the following statement is false for
$\varepsilon\in\{0,1\}$\thinspace:

\smallskip {\bf Theorem 1.1b.} {\it Let $(X,d)$ be a complete metric space and
$f:X\rightarrow X$ a self mapping of $X$. Let $\Lambda\geq0$,
$\alpha\geq1$ and $\beta\in\lbrack0,\alpha]$ be fixed constants. If
there exists $\varepsilon\in\lbrack0,1]$ such that
\begin{equation}
\forall x,y\in X:d(f(x),f(y))\leq(1-\varepsilon)d(x,y)+\Lambda\varepsilon
^{\alpha}\psi(\varepsilon)[1+\left\Vert x\right\Vert +\left\Vert y\right\Vert
],\label{r-b1.2b}%
\end{equation}
then $f$ possesses a (unique) fixed point $x^{\ast}=f(x^{\ast})$.}

\medskip
\medskip In fact Theorem 1.1b is false not only for $\varepsilon\in\{0,1\};$
it is false for any fixed $\varepsilon\in\lbrack0,1]$. Indeed, take,
for example, $X:=\{-1,1\}\subset(\mathbb{R},\left\vert
\cdot\right\vert )$ with $f(x):=-x$. Fix some $\alpha\geq1$,
$\psi:[0,1]\rightarrow\mathbb{R}_{+}$ (increasing, continuous at $0$
with $\psi(0)=0$), and $\varepsilon\in \lbrack0,1]$. Clearly
$\left\vert f(x)-f(y)\right\vert =\left\vert x-y\right\vert
\leq(1-\varepsilon)d(x,y)+\varepsilon\left(  1+\left\Vert
x\right\Vert +\left\Vert y\right\Vert \right)  $ for $x,y\in X$.
Setting $\Lambda:=\varepsilon^{\alpha-1}\psi(\varepsilon)$ for
$\varepsilon\in(0,1]$ and $\Lambda=0$ for $\varepsilon=0$,
(\ref{r-b1.2b}) is verified. Clearly, $f$ has not fixed points.

\medskip Recall also that besides the main aim of \cite{Ber18}, that is to
prove that \textquotedblleft Theorem 1.1 does not hold if we have
$\varepsilon=0$ or $\varepsilon=1$ in (1.2)", another aim was to
\textquotedblleft provide a correct version (but not fully in the
spirit) of Theorem 1.1".

\smallskip That \textquotedblleft correct version (but not fully in the
spirit) of Theorem 1.1" is \textquotedblleft the following interesting
existence result.

\smallskip {\bf Theorem 3.5.} {\it Let $(X,d)$ be a complete metric space and
$f:X\rightarrow X$ a self mapping of $X$. Let $\Lambda\geq0$,
$\alpha\geq1$ and $\beta\in\lbrack0,\alpha]$ be fixed constants. If
the inequality

$d(f(x),f(y))\leq(1-\varepsilon)d(x,y)+\Lambda\varepsilon^{\alpha}%
\Psi(\varepsilon)d(y,f(x))\ \ \ (3.16)$

\noindent is satisfied for some $\varepsilon\in(0,1)$ and every
$x,y\in X$, then

$1)$ $Fix(T)=\{x\in X:Tx=x\}\neq\emptyset;$

$2)$ For any $x_{0}\in X$, Picard iteration
$\{x_{n}\}_{n=0}^{\infty}$, $x_{n}=T^{n}x_{0}$, converges to some
$x^{\ast}\in FixT$;

$3)$ The following estimate holds

$d(x_{n+i-1},x^{\ast})\leq\frac{(1-\varepsilon)^{i}}{\varepsilon}%
d(x_{n},x_{n-1})$ \ \ $n=0,1,2...;$ $i=1,2,...\ \ \ (3.17)$}".

\medskip

Related to \cite[Theorem 3.5]{Ber18} let us observe the following:
1) Surely, one must have $f$ instead of $T$. 2) The introduction of
$\beta$ in its statement is superfluous because $\beta$ is neither
involved in Eq.\ (3.16), nor in its conclusion; even more, the
introduction of $\alpha$ is superfluous because for
$\varepsilon\in(0,1)$ fixed, setting $\delta:=1-\varepsilon$ and
$L:=\Lambda\varepsilon^{\alpha}\psi(\varepsilon)$, the hypothesis of
\cite[Theorem 3.5]{Ber18} is saying that $f$ is \textquotedblleft an
almost contraction" (see \cite[Definition 3.1]{Ber18}), and so
\cite[Theorem 3.5]{Ber18} is nothing but \cite[Theorem 3.4]{Ber18}.

\medskip As a conclusion, one can say that \cite[Theorem 1.1]{Ber18} has only
a visual similarity with \cite[Theorem 1]{Pat11}. It seems that the aim of
\cite{Ber18} is to advertise its author's results on fixed point theorems
because we may not admit he does not know the meaning of the quantifiers
$\forall$ and $\exists$.

\medskip

We end this note with the following remark extracted from \cite{ChKaMeRa19}:

\begin{quotation}
\textquotedblleft Remark 1 There is a scope of misunderstanding with the
inequality (1) in the work of Pata and also in similar other inequalities like
(2)--(4) in works incorporating the ideas of Pata. Berinde noted in [7] that
if the condition (1) is satisfied, not for all $\varepsilon\in\lbrack0,1]$,
but just for some specific values, the conclusion of Theorem 1 might not hold.
For example, if (1) holds just for $\varepsilon=0$, then one has just the
non-expansive condition

$d(fx,fy)\leq d(x,y)\ \ \forall x,y\in X$, \noindent which obviously
does not imply the existence of fixed point. Similarly, if
$\varepsilon=1$, (1) reduces to

$d(fx,fy)\leq L[1+\left\Vert x\right\Vert +\left\Vert y\right\Vert
]^{\beta}$, \noindent with some constant $L$, which is also known to
be insufficient for the existence of a fixed point of $f$.

Similar conclusions hold for conditions (2)--(4).

From this observation Berinde concludes that the Pata-type result is
incorrect. But this is not so. There is no contradiction between the above
observations and conclusions of the Pata-type theorems for the following
reasons. As we have already noted in the introduction, the Pata type results
are obtained for functions satisfying a family of inequalities and any single
inequality from the above mentioned family will not provide us with a
sufficient condition for the existence of a fixed point. Had it been so, then
there is no need of considering a family of inequalities. Thus the argument of
Berinde is not tenable."
\end{quotation}

\end{document}